\documentclass[12pt]{article}
\usepackage{amsxtra,amssymb,amsthm,amsmath,latexsym}

\theoremstyle{plain}
\newtheorem{theorem}{Theorem}

\def\R{{\mathbb R}}

\def\oH{{\overset{\circ}{H}}}
\def\oH1{{\overset{\circ}{H}\kern-.02in{}^1}}

\def\bee{\begin{equation*}}
\def\eee{\end{equation*}}
\def\be{\begin{equation}}
\def\ee{\end{equation}}

\begin{document}

\title{Inverse problems for parabolic equations 3 }

\author{A.G. Ramm\\
 Mathematics Department, Kansas State University, \\
 Manhattan, KS 66506-2602, USA\\
ramm@math.ksu.edu }

\date{}
\maketitle\thispagestyle{empty}

\begin{abstract}

Let $u_t-a(t)u_{xx}=f(x, t)$ in $0\leq x \leq \pi,\,\,t\geq 0.$
Assume that $u(0,t)=u_1(t)$,  $u(\pi,t)=u_2(t)$,  $u(x,0)=h(x)$,
and the extra data $u_x(0,t)=g(t)$ are known.
The inverse problem is: {\it

How does one determine the unknown $a(t)$?}

The function $a(t)>a_0>0$ is assumed continuous and bounded.
This question is answered and a method for recovery
of $a(t)$ is proposed.  There are several papers
in which sufficient conditions are given for
the uniqueness and existence of $a(t)$,  but apparently there was no
method  proposed for calculating of $a$.
The method given in this paper for proving the uniqueness and existence of the solution
to inverse problem is new and it allows one to calculate the unknown coefficient $a(t)$.
\end{abstract}

\section{Introduction}\label{S:1}

Consider the problem
\be\label{e1}u_t-a(t)u_{xx}=f(x,t)
\quad (x,t)\in
[0,\pi]\times[0, \,\infty), \ee
\be\label{e2}u(0,t)=u_1(t), \quad u(\pi, t)=u_2(t), \quad u(x,0)=h(x),
\ee
where the functions $\{u_1,\,u_2,\, h\}$ are known. The extra data are
\be\label{e3}
u_x|_{x=0}=g(t).
\ee
The Inverse Problem (IP) we are interested in is the following one:

{\it How can one find $a(t)$ given the data?}

There is an extensive literature on inverse problems
for the heat equation (see \cite{I}, \cite{R1},  and references therein),
but there was no method for calculating $a(t)$  as far as the author knows.
The method given in this paper for proving the uniqueness and existence of the solution
to inverse problem is new and it allows one to calculate the unknown coefficient $a(t)$.
In \cite{R2}-\cite{R7} the author studied various inverse problems for
parabolic equations.

Let $||u||:=||u||_{L^2(0,\pi)},\,\, u_m:=(u, p_m)=\int_0^{\pi}up_m dx$,
\be\label{e4}
p_m^{\prime\prime}+m^2p_m=0,\quad 0\leq x \leq \pi,\quad
p_m(0)=p_m(\pi)=0,\quad
||p_m||=1,
\quad m=1,2\dots,
\ee
where $p_m=\sqrt{\frac 2 {\pi}}\sin (mx)$.
Let
\be\label{e5}
u=v+u_1+\frac{x}{\pi}(u_2-u_1):=v+r.
\ee
Then $v$ solves the problem:
\be\label{e6}
v_t-a(t)v_{xx}=f-r_t:=F(x,t), \quad v(0,t)=v(\pi,t)=0, \quad v(x,0)=h(x)-r(x,0):=H(x).
\ee
Let us look for the solution to problem (6) of the form
\be\label{e7}
v(x,t)=\sum_{m=1}^\infty c_m(t)p_m(x),
\ee
where the functions $c_m(t)$ are to be found. Define
\be\label{e8}
H(x)=\sum_{m=1}^\infty H_m p_m(x),
\ee
and
\be\label{e9}
F(x,t)=\sum_{m=1}^\infty F_m(t) p_m(x),
\ee
Multiplying equation (6) by $p_m(x)$ and integrating
over the interval $[0,\pi]$ and then by parts, using the boundary conditions (6), one gets
\be\label{e10}
\dot {c}_m+a(t)m^2c_m=F_m h(t),\quad c_m(0)=H_m,\quad
\ee
where $m=1,2, \dots$.
Define
\be\label{e11}
A(t)=\int_0^ta(s)ds.
\ee
 The unique solution to problem (10) is:
 \be\label{e12}
 c_m(t)=H_me^{-m^2 A(t)}+e^{-m^2 A(t)}\int_0^t e^{m^2 A(s)}F_m(s)ds.
 \ee
Therefore $v$ and $u$ can be found:
\be\label{e13}
u=\sum_{m=1}^\infty \Big( H_me^{-m^2 A(t)}+e^{-m^2 A(t)}\int_0^t e^{m^2 A(s)}F_m(s)ds\Big)p_m(x)+r(x,t).
\ee
Let us use extra data (3):
\be\label{e14}
\sum_{m=1}^\infty \Big( H_me^{-m^2 A(t)}+e^{-m^2 A(t)}\int_0^t e^{m^2 A(s)}F_m(s)ds\Big)m(2/\pi)^{1/2} +[u_2-u_1]/T =g(t).
\ee
Define two functions:
\be\label{e15}
Q_0(z):=\sum_{m=1}^\infty m H_me^{-m^2 z},
\ee
and
\be\label{e16}
Q(z,s):=\sum_{m=1}^\infty m F_m(s)e^{-m^2 z}.
\ee
Then equation (14) takes the form:
\be\label{e17}
Q_0(A(t))+\int_0^tQ(A(t)-A(s),s)ds=g(t).
\ee
This is an equation for $A(t)$, and $a(t)=dA(t)/dt$. Assume that
\be\label{e18}
H_m>0   \quad \forall m\ge 1.
\ee
Then $dQ_0(z)/dz>0$, so the function $Q_0$ has an inverse function $Q^{-1}_0: \R_+\to R_+,$ where $\R_+=(0,\infty)$.
Inverting $Q_0$ in equation (17) one gets:
\be\label{e19}
A(t)=-Q^{-1}_0\int_0^tQ(A(t)-A(s),s)ds+Q^{-1}_0g(t):=V(A)+G:=T(A).
\ee
Here we assume that
\be\label{e20}
g(t)> 0,
\ee
and
\be\label{e21}
F_m(t)\ge 0   \quad \forall m\ge 1,
\ee
which is sufficient for the inequality $Q(A(t)-A(s),s)>0$ to hold.
Equation (19) is a Volterra-type equation for positive $A(t)$ and the integral operator
in (19) is an operator which maps a cone of non-negative functions into itself.
One has
\be\label{e22}
dQ/dz<0, \quad max_{z,s}|\frac{\partial Q(z,s}{\partial z}|\le C,\quad C=const,
\ee
provided that
\be\label{e23}
max_{t\ge 0} \sum_{m=1}^\infty m^3|F_m(t)|\le C' \quad C'=const.
\ee
  Our result is the following theorem.

\begin{theorem}\label{T:1}
Assume (18), (20), (21) and (23). Then $A(t)>0$ is uniquely determined  by equation
(19). This equation can be uniquely solved by iterations
\be\label{e24}
A_{n+1}(t)=-Q^{-1}_0\int_0^tQ(A_n(t)-A_n(s),s)ds+Q^{-1}_0g(t), \quad A_1(t)=Q^{-1}_0g(t).
\ee
 \end{theorem}

In Section 2 proof is given.

\section{Proofs}\label{S:2}

\begin{proof}[Proof of Theorem 1.]
Only convergence of the process (24) is not yet proved. This convergence locally, for small $t$, $t\le t_0$,
 where $t_0$ is a sufficiently small number, follows from the
contraction mapping principle. To apply this principle one chooses a ball $B_M$ in the space of continuous
functions with the usual norm,
$max_{s\in [0,t_0]} |A(s)|\le M$, uses the estimates $0<Q_0^{-1}(z)\le C_0$, $0<Q(z,s)\le C$, $|\frac{\partial Q}{\partial z}|\le C_1$, where $C_0, C, C_1$ are some constants, and checks that for sufficiently small $t$ the operator $T$ in equation (19) is a contraction on the ball $B_M$. Global convergence ($t_0=\infty$)  follows
from the argument  usual for Volterra-type equations. Although the operator $V$ in equation (19) is a
nonlinear Volterra-type operator, but since $C_1$ is finite, the usual estimates hold. This
proves, under our assumptions, global existence and uniqueness of the positive solution to equation (19).
The solution $A$ to this equation is not only positive but monotone increasing. Indeed, the function
$Q_0$ is monotone decreasing, so $Q_0^{-1}$ is monotone increasing. The integral in (19) is monotone
increasing since the integrand is non-negative. Therefore, $A(t)$ is monotone increasing. Consequently,
$a(t)=dA(t)/dt \ge 0$.

Theorem 1 is proved. \hfill $\Box$

\end{proof}

{\bf Remark 1.} Let us point out some cases when the inverse problem can be solved explicitly,
in closed form. Assume that $u_1=u_2=h=0$ and $f(x,t)=p_1(x)$.
Note that the choice of $u_1, u_2, h$ and $f$ can be made by the experimenter so that the unknown
coefficient $a(t)$ can be found easier.
Look for the solution $u=c(t)p_1(x)$.
Then $c_t+a(t)c=1$, $c(0)=0$, so $a(t)=(1-c_t)/c$. The extra data are $g(t):=u_x|_{x=0}=c(t)(2/\pi)^{1/2}$.
Therefore $c(t)=g(t)(\pi/2)^{1/2}$, so $a(t)=((2/\pi)^{1/2}-g(t))/g'(t)$. Other cases when the unknown
functions can be expressed in closed form through the data  can be given.
$$ $$

{\it Acknowledgement.} The author thanks Prof. M.Ivanchov
who drew his attention to this inverse problem.


\begin{thebibliography}{1000} 

 \bibitem{I} Ivanchov, M. I., {\bf Inverse Problems for Equations of Parabolic Type}, VNTL, Lviv, 2003.

\bibitem{R1} 
Ramm, A.~G.,
{\bf Inverse Problems}, Springer, New York, 2005.

\bibitem{R2}
Ramm, A.~G., {\it An inverse problem for the heat equation},
J.~Math.~Anal.~Appl., 264, N2, (2001), 691-697.

\bibitem{R3}
Ramm, A.~G., {\it Inverse problems for parabolic equations},
Australian Jour. ~Math.~Anal.~Appl., 2, N2, (2005).

\bibitem{R4}
Ramm, A.~G., {\it An inverse problem for the heat equation},
 Proc.Roy.Soc. Edinburgh, 123, N6, (1993), 973-976.

\bibitem{R5}
Ramm, A.~G., Inverse problems for parabolic equations 2,
 Communic. in Nonlinear Sci. and Numer. Simulation,
12, (2007), 865-868.

\bibitem{R6} Ramm, A.~G., Property $C$ for ODE  and applications to an
inverse problem for a heat equation,
Bull. Polish Acad of Sci. Mathem., 57, N3-4, (2009), 243-249.

\bibitem{R7} Ramm, A.~G., An inverse problem for the heat equation II,
Applic. Analysis, 81, N4, (2002), 929-937.


\end{thebibliography}
\end{document}